\documentclass[12pt,leqno]{article}
\usepackage{amstext,amsmath,amsthm,amsfonts,amssymb,amscd}
\usepackage{graphics}
\usepackage[latin1]{inputenc}

\usepackage{graphicx}
\usepackage{color}

\newcommand\bib[1]{\bibitem[#1]{#1}}
\begin{document}

\centerline{\LARGE\bf Parallel Mean Curvature Surfaces in }

\medskip

\centerline{\LARGE\bf Symmetric Spaces}

\vskip .4in

\centerline{\bf by}

\vskip .2in

\centerline{\large\bf  M. J. Ferreira\, and \, Renato Tribuzy}

\vskip .5in

\begin{abstract}
We present a reduction of codimension theorem for surfaces with
parallel mean curvature in symmetric spaces.
\end{abstract}

\footnotetext{MSC 2000 -- primary: 53C42. 53C40;
 secondary: 53A10}

\noindent{\large\bf I - Introduction}

\medskip

Surfaces with parallel mean curvature have been considered by many
geometers. Yau [Y] studied them when the ambient space is a space
form. He showed that such a surface must be, either minimal in a
round sphere, or its image must be contained in a three dimensional
totally geodesic, or totally umbilic submanifold. This result was
recently extended to maps into $E^n(c) \times \mathbb R$ by Alencar,
do Carmo and Tribuzy [A-C-T], where $E^n(c)$ a space form with
constant sectional curvature $c \neq 0$. They have proved that the
image of the surface is either minimal in a totally umbilical
hypersurface of $E^n(c)$, has constant mean curvature in a three
dimensional totally geodesic, or totally umbilical, submanifold of
$E^n(c)$, or lies in $E^4 \times \mathbb R$.

Fetcu [F] has shown that a surface immersed in $CP^n$ with parallel
mean curvature must be either pseudo-umbilic and totally real or its
image must be contained in $CP^5$. We recall that a immersion
$\varphi: M \rightarrow N$ between two Riemannian manifolds is said
to be pseudo-umbilic if the $H$-Weigarten operator $A_H$ is a
multiple of the identity, where $H$ denotes the mean curvature.

The aim of this work is to generalize the above results to
immersions into symmetric spaces.

From now on we let $S$ denote a Riemannian symmetric space and $R$
its curvature tensor. For $p \in S$, we say that a subspace $K
\subset T_pS$ is invariant by the curvature tensor $R_p$, if
$R_p(u,v)w \in K$, whenever $u,v,w \in K$. When $X \subset T_p S$,
$R_p(X)$ will represent the least vector subspace of $T_p S$
invariant by the curvature tensor $R_p$ at $p$.

Let $\varphi : M \rightarrow S$ be an isometric immersion from a
Riemann surface $M$. For each $x \in M$, the $n$-th osculating space
of $\varphi$ at $x$ will be denoted by $O^n_x (\varphi)$.

Our main result is the following:

\bigskip

\noindent{\bf Theorem 1}

If $\varphi$ has parallel mean curvature $H\neq 0$, one of the
following conditions hold:

i) $\varphi$ is pseudo-umbilic;

ii) The dimension $d$ of $R_{\varphi(x)}(O^2_{x})$ is independent of
$x \in M$ and there exists a totally geodesic $S^{\prime} \subset
S$, with dimension $d$, in such a way that $\varphi (M) \subset
S^{\prime}$ .

\medskip

\bigskip

\noindent{\large\bf II - Proof of the main theorem}

\medskip

The following theorem dealing with reduction of codimension of maps
is fundamental to our purpose:

\noindent{\bf Theorem 2}

Let $\varphi:M \rightarrow S$ be an isometric immersion from a
Riemannian manifold into a symmetric space. If there exists a
parallel fiber bundle $L$, over $M$, such that $R(L)=L$ and $TM
\subset L$, then there exists a totally geodesic submanifold
$S^{\prime}$ of $S$ with $\varphi (M) \subset S^{\prime}$ and $L_x =
T_{\varphi(x)}S^{\prime}$, $\forall x \in M$.

\medskip

For a proof see [E-T].

\bigskip

Let $\alpha$ denote the second fundamental form of the immersion. We
remark that, for each $x \in M$, $O^2_x(\varphi)=T_x M + N_1(x)$,
where $N_1(x) =\{\alpha(u,v): u,v \in T_x M\}$. We will show that
either $H$ is a umbilical direction, or $L=R(O^2(\varphi))$ is a
parallel fiber bundle.

Let $U$ be an open subset of $M$ where $L_x = R_x (O^2_x (\varphi))$
has maximal dimension. We will prove that $L$ is parallel on $U$
when $H$ is not a umbilical direction. Using analiticity arguments
we may conclude that its dimension is constant on $M$.

\bigskip

\noindent{\bf Lemma 1}

Whenever $X,Y,W \in \Gamma(L)$ and $\nabla X, \nabla Y, \nabla W \in
\Gamma (L)$, $R(X,Y)W \in \Gamma (L)$.

\bigskip

\noindent{\bf Proof of Lemma 1:}

It follows straightforward from the parallelism of $R$.

\bigskip

We remark that $L$ is constructed from $O^2(\varphi)$ by successive
applications of $R$ to its elements. Thus it is enough to show that
$ \nabla \xi \in \Gamma(L)$, for every $\xi= \alpha (U,V)$, to
conclude that $L$ is parallel.

Let $\{\varepsilon_1, \varepsilon_2\}$ be a local orthonormal frame
field defined on open subset $V$ of $U$. Notice that
$R_x(\varepsilon_1,\varepsilon_2)H \in O^2_x$, $\forall x\in V$. Let
$\pi_{N_1(x)}:T_{\varphi(x)} S \rightarrow N_1(x)$ denote the
orthogonal projection and consider
$n_1(x)=\pi_{N_1(x)}(R_x(\varepsilon_1,\varepsilon_2)H)$.

\bigskip

\noindent{\bf Lemma 2}

Assume $n_1$ vanishes identically on $V$. Then $L$ is parallen on
$V$.

\bigskip
\noindent{\bf Proof of Lemma 2:}

From Ricci equations we get

$$
0=<R^{\perp}(\varepsilon_1,\varepsilon_2)H,\eta>=-<[A_H,
A_\eta]\varepsilon_1,\varepsilon_2>,
$$
so that $A_H$ commutes with $A_\eta$ for every section $\eta$ of the
normal bundle. Therefore without loss of generality
$\{\varepsilon_1,\varepsilon_2\}$ may be chosen in such a way that
$\alpha(\varepsilon_1,\varepsilon_2)=0$ and it is enough to show
that
$$
\nabla_{\varepsilon_i} \alpha(\varepsilon_j,\varepsilon_j)\in
\Gamma(L)
$$.

But this is equivalent to prove that

$$
(\nabla_{\varepsilon_i} \alpha)(\varepsilon_j , \varepsilon_j) \in
\Gamma(L),
$$
since

$$
(\nabla_{\varepsilon_i} \alpha)(\varepsilon_j ,
\varepsilon_j)=\nabla_{\varepsilon_i}
\alpha(\varepsilon_j,\varepsilon_j)- 2
\alpha(\nabla_{\varepsilon_i}\varepsilon_j,\varepsilon_j).
$$

From Codazzi equations we have:

$$
(\nabla_{\varepsilon_i} \alpha)(\varepsilon_j ,
\varepsilon_j)=(\nabla_{\varepsilon_j} \alpha)(\varepsilon_i ,
\varepsilon_j)+(R(\varepsilon_i,\varepsilon_j)\varepsilon_j)^{\perp},
$$
where $( ~~  )^{\perp}$ represents the orthogonal projection onto
the normal bundle. But
$(R(\varepsilon_i,\varepsilon_j)\varepsilon_j)^{\perp}$ sits in $L$
since
$$
(R(\varepsilon_i,\varepsilon_j)\varepsilon_j)^{\perp}=R(\varepsilon_i,\varepsilon_j)\varepsilon_j
- A,
$$
with $A=(R(\varepsilon_i,\varepsilon_j)\varepsilon_j)^{T} \in
\Gamma(L)$, where $(~~ )^T$ stands for the orthogonal projection
onto $TV$.

Of course $\nabla_{\varepsilon_i}\alpha(\varepsilon_j,\varepsilon_j)
\in \Gamma(L)$ if $i\neq j$. When $i=j$,
$$
\nabla_{\varepsilon_i} \alpha(\varepsilon_i,\varepsilon_i)=
\nabla_{\varepsilon_i}(2H-\alpha(\varepsilon_j,\varepsilon_j))=-\nabla_{\varepsilon_i}
\alpha(\varepsilon_j,\varepsilon_j)+B \in \Gamma(L),
$$
with $j\neq i$ and $B=2(\nabla_{\varepsilon_i}H)^T$.
\bigskip

We proceed now with the proof of theorem 1, considering the case
where $n_1$ does not vanish identically.

Let us consider an open subset $V$ of $U$ where $n_1 (x) \neq 0$,
$\forall x \in V$.

We remark that $n_1$ and $H$ are linearly independent, $\nabla n_1
\in \Gamma(L)$ and $\nabla H \in \Gamma(L)$. Remember that, due to
the parallelism of $R$, $\nabla n_1=R(\nabla
\varepsilon_1,\varepsilon_2)H+R(\varepsilon_1,\nabla \varepsilon_2)H
+ R(\varepsilon_1,\varepsilon_2)\nabla H$.

Let $Z$ denote the sub-bundle of the normal bundle spanned by $n_1$
and $H$. Consider a unitary section $\gamma$ orthogonal to $Z$ and
sitting in $N_1$ and take an orthonormal frame field
$\{\varepsilon_1,\varepsilon_2\}$ diagonalizing the Weingarten
operator $A_{\gamma}$.

We observe that $\alpha(X,Y)=\alpha(X,Y)_{Z}+\alpha(X,Y)_{\gamma}$,
where $\alpha_Z$ and $\alpha_\gamma$ are respectively the $Z$ and
$<\gamma>$ components of $\alpha$; here $<\gamma>$ denotes the
sub-bundle spanned by $\gamma$.

Now, to show that $\nabla \alpha(X,Y) \in \Gamma(L)$, it is enough
to prove that $\nabla \alpha_\gamma(X,Y) \in \Gamma(L)$, since
$\alpha_Z(X,Y) \in \Gamma(L)$. Therefore we get that, for any $i,j,k
\in \{1,2\}$, $\nabla_{\varepsilon_i} \alpha_\gamma(\varepsilon_j,
\varepsilon_k) \in \Gamma(L)$ if and only if $\nabla_{\varepsilon_i}
\alpha(\varepsilon_j, \varepsilon_k)$, which is equivalent to
$(\nabla_{\varepsilon_i} \alpha)(\varepsilon_j, \varepsilon_k)\in
\Gamma(L)$.

Clearly
$\nabla_{\varepsilon_i}\alpha_\gamma(\varepsilon_j,\varepsilon_k)=0$,
whenever $j\neq k$.

When $j=k \neq i$, using Codazzi equations, we have
$$
(\nabla_{\varepsilon_i}\alpha)(\varepsilon_j,\varepsilon_j)=(\nabla_{\varepsilon_j}\alpha)(\varepsilon_i,\varepsilon_j)+A,
$$
where $A\in \Gamma(L)$ and, since
$(\nabla_{\varepsilon_j}\alpha_{\gamma})(\varepsilon_i,\varepsilon_j)
\in \Gamma(L)$, we obtain
$(\nabla_{\varepsilon_i}\alpha)(\varepsilon_j,\varepsilon_j)\in
\Gamma(L)$.

If $i=j=k$,

$$
\nabla_{\varepsilon_i}\alpha_\gamma(\varepsilon_i,\varepsilon_i)=\nabla_{\varepsilon_i}(2H_\gamma
-
\alpha_\gamma(\varepsilon_j,\varepsilon_j))=-\nabla_{\varepsilon_i}\alpha_\gamma(\varepsilon_j\varepsilon_j),
$$
where $j \neq i$ and we have used the fact that $H_\gamma = 0$.
Hence
$\nabla_{\varepsilon_i}\alpha_\gamma(\varepsilon_i,\varepsilon_i)
\in \Gamma(L)$.

\medskip

\bigskip

\noindent{\large\bf III - Remarks}

\medskip

\noindent{\bf 1 -}

When $S=E^n(c)$ is a space form of constant sectional curvature $c$,
it follows directly from Ricci equation that $A_H$ is either a
multiple of the identity, or there exists a basis diagonalizing, at
each point, the second fundamental form. Hence the second normal
space has dimension less or equal than $4$. Of course the curvature
tensor leaves this space invariant and we get that $\varphi(M)
\subset E^4(c)$ ([Y]).

\bigskip

\noindent{\bf 2 -}

When $S=CP^n$ we observe that
$$
R(O^2(\varphi)) \subset O^2(\varphi) + J O^2(\varphi),
$$
where $J$ denotes the standard complex structure of $CP^n$.
Therefore, if $\varphi$ is not pseudo-umbilic, we must clearly have
$\varphi(M) \subset CP^5$. In the case were $H$ is a umbilical
direction, a direct computation gives $R(X,Y)H=0$, whenever $X,Y \in
\Gamma(TM)$. Thus, from the explicit formula of the curvature
tensor, we obtain straightforward that $JTM$ is orthogonal to $TM$,
i.e $\varphi$ is totally real. This case was studied by Fetcu [F].

\bigskip

\noindent{\bf 3 -}

Assume now that the target manifold is a product $S_1 \times S_2$ of
symmetric spaces. We have:

\bigskip
\noindent{\bf Corolary 1 - }
Let $\varphi :M \rightarrow S_1 \times
S_2$ an immersion from a Riemann surface with parallel mean
curvature $H \neq 0$. Then, one of the following conditions hold:

i) $\varphi$ is pseudo-umbilic;

ii) $\varphi(M) \subset S_1^{\prime} \times S_2^{\prime}$, where,
for each $i\in \{1,2\}$, $S_i^{\prime}\subset S_i$ is a symmetric
space totally geodesic embedded in $S_i$ and $\dim S_i^{\prime}=dim
R_i((O^2(\varphi))^{T_i})$; here $R_i$ represents the curvature
tensor of $S_i$ and $(~~)^{T_i}$ stands for the orthogonal
projection onto $TS_i$.

\medskip

\noindent{\bf Proof of Corollary 1:}

\medskip

Clearly
$$
R(O^2(\varphi) )\subset R_1((O^2(\varphi))^{TS_1}) +
R_2((O^2(\varphi))^{TS_2}).
$$
Therefore it is enough to prove that each
$R_i((O^2(\varphi))^{TS_i})$ is parallel. Let us take an open subset
of $M$ where the dimension of $O^2(\varphi)^{T_i}$ is maximal. By
Lemma 1, it suffices to show that $\nabla \mu_i \in \Gamma
(R((O^2(\varphi))^{T_i}))$, whenever $\mu_i \in
\Gamma(O^2(\varphi))$. This follows from the fact that there exists
$\mu \in\Gamma(O^2(\varphi))$ such that $\mu_i = (\mu)^{T_i}$. Since
$\nabla \mu = (\nabla \mu_1, \nabla \mu_2) \in R(O^2(\varphi))$ we
obtain $\nabla \mu_i \in R_i(O^2(\varphi))$.

\bigskip

In the particular situation where $S_2=\mathbb R$, it follows from
corollary 1 that, either $\varphi$ is pseudo-umbilic, or $\varphi(M)
\subset S^{\prime} \times \mathbb R$, where $\dim S^{\prime}=\dim
R((O^2(\varphi))^ {TS})$ and $R$ denotes the curvature tensor of
$S$. For instance, if $S=E^n(c)$, clearly $S^{\prime}=E^4(c)$, like
in [A-C-T].

When the target manifold is $E^n(c) \times E^m(d)$, either $\varphi$
is pseudo-umbilic, or $\varphi(M) \subset E^4(c) \times E^4(d)$.


\bigskip

mjferr@ptmat.fc.ul.pt

tribuzy@pq.cnpq.br

\bigskip

Work supported by FCT (Funda\c{c}\~{a}o para a Ci\^{e}ncia
 e Tecnologia) in Portugal and by CNPq and FAPEAM in Brasil

\end{document}